\documentclass[12pt]{article}
\usepackage{amsfonts,amssymb,amsmath}
\usepackage{graphicx,caption}

\def\R{\mathbb{R}}

\def\p{ \partial }

\def\bq{ \begin{equation} }
\def\eq{ \end{equation} }
\def\ben{ \begin{eqnarray} }
\def\en{ \end{eqnarray} }

\def\frac#1#2{{#1\over #2}}
\def\on#1#2{\mathop{\vbox{\ialign{##\crcr\noalign{\kern2pt}
$\scriptstyle{#2}$\crcr\noalign{\kern2pt\nointerlineskip}
\kern-2pt$\hfil\displaystyle{#1}\hfil$\crcr}}}\limits}

\begin{document}

\title{Linearisable Abel equations \\ and the Gurevich--Pitaevskii problem}
\author{S. Opanasenko$^1$, E.V. Ferapontov$^{1, 2}$}
   \date{}
\vspace{-20mm}
   \maketitle
\vspace{-7mm}
\begin{center}
${^1}$ Department of Mathematical Sciences \\
Loughborough University \\
Loughborough, Leicestershire LE11 3TU, UK \\[1ex]
    \ \\
$^2$Institute of Mathematics, Ufa Federal Research Centre\\
Russian Academy of Sciences, 112 Chernyshevsky Street \\
Ufa 450008, Russian Federation\\
e-mails: \\
\texttt{S.Opanasenko@lboro.ac.uk}\\
\texttt{E.V.Ferapontov@lboro.ac.uk}
\end{center}

\medskip

\begin{abstract}
Applying symmetry reduction to a class of $\mathrm{SL}(2,\mathbb R)$-invariant third-order ODEs, we obtain Abel equations
whose general solution can be parametrised by hypergeometric functions.
Particular case of this construction provides a general parametric solution to the Kudashev equation, an ODE arising in the Gurevich--Pitaevskii problem, thus giving  the first term of a large-time asymptotic expansion of its solution in the oscillatory (Whitham) zone.

\medskip


\medskip

\end{abstract}

\newpage

\section{Introduction}
The first-kind and the second-kind Abel equations,
\[
u'=f_3(x)u^3+f_2(x)u^2+f_1(x)u+f_0(x)
\]
and
\[
(G_1(x)U+G_0(x))U'=F_3(x)U^3+F_2(x)U^2+F_1(x)U+F_0(x),
\]
which are related by a point transformation $u=1/(G_1U+G_0)$, have been of interest since the work of  Abel on elliptic functions
where Abel equations appeared in the form
\[
(\eta+h_0(x))\eta'=H_2(x)\eta^2+H_1(x)\eta+H_0(x).
\]
Although Abel equations can be seen as a slight generalization of well-understood Riccati equations,
there is no general approach to solving them.
By and large, there has been no real progress on this issue since the classical works of Abel, Liouville and Appell
at the end of the 19th century \cite{Abel, Appell, Liouville}.
In fact, the classical textbooks  on solutions of ODEs (Kamke \cite{Kamke1959}, Polyanin and Zaitsev \cite{PolyaninZaitsev2003}) contain about a hundred of integrable Abel equations,
however, most of them are equivalent to the 11 canonical  forms \cite[Appendix~A]{ChebTerrabRoche2000} under the transformations
$$
x\to \phi(x), \quad u\to \psi(x)u+\eta(x);
$$
note that under these transformations every Abel equation of the first kind can be reduced to the normal form $u'=u^3+f(x)$. We also note that Chiellini integrability condition, a discovery of 1930s and a working horse of some modern progress~\cite{MakHarko2002,MancasRosu2013},
is merely a means of checking whether an Abel equation at hand is equivalent to a separable Abel equation.

In the present paper we obtain a two-parameter family
 of integrable Abel equations as a reduction of $\mathrm{SL}(2,\mathbb R)$-invariant third-order ODEs (examples of solvable  Abel equations obtained as symmetry reductions of some second-order ODEs can be found in \cite{Mustafa}). Although the found equations possess first integrals which are clear generalizations of that of known integrable Abel equations,
these integrals are quite cumbersome, and we present more palatable parametric solutions in terms of hypergeometric functions.
As an application, we provide a general parametric solution to the Kudashev equation,
\begin{gather}\label{eq:Kudashev}
\frac{\mathrm dR}{\mathrm dz}=\frac{486R^4-171R^2+9zR+5}{9(54R^3-9R+z)(2R+3z)},
\end{gather}
which takes a form of a second-kind Abel equation when written in terms of $z(R)$.
Equation~(\ref{eq:Kudashev}) arises in the following context (see \cite{GST, SS} and references therein): consider the KdV equation
constrained by the stationary part of its higher-order non-autonomous symmetry,
\begin{equation}\label{kdv}
u_t+uu_x+u_{xxx}=0, \qquad u_{xxxx}+\frac{5}{3}uu_{xx}+\frac{5}{6}u_x^2+\frac{5}{18}(x-tu+u^3)=0.
\end{equation}
Note that the second equation (\ref{kdv}),  denoted $P_I^2$, is often viewed as a fourth-order analogue of the classical Painlev\'e I equation.
Under the Gurevich-Pitaevskii boundary condition \cite{GP}, namely, $u \to -x^{1/3}$ as $x\to \infty$,
 equations (\ref{kdv})  possess a unique real solution with no poles on the real line \cite{Claeys}.
 This solution  possesses several different asymptotic expansions depending on the value of the self-similar variable $z=xt^{-\frac{3}{2}}$, see \cite{Claeys10} for the precise statements and the description of the Riemann-Hilbert problem associated with the real pole-free solution of the equation~$P_I^2$, which was analysed based on the Deift/Zhou steepest descent method \cite{DZ}. The large-time
 asymptotic expansion in the oscillatory zone was previously discussed in \cite{Potemin} based on the Whitham averaging theory. An alternative approach, which does not refer to the Whitham averaging theory,  was developed in \cite{GST} where one seeks a large-time asymptotic expansion in the form
\begin{equation}\label{asy}
u(t,x)=\sqrt t\left(v_0(z, \phi)+t^{-\frac{7}{4}}v_1(z, \phi)+t^{-\frac{7}{2}}v_2(z, \phi)+ \dots \right).
\end{equation}
Here $z=xt^{-\frac{3}{2}}$ and $\phi=t^{\frac{7}{4}}f(z)+S(z)$ are the slow and the fast variables, respectively
(the functions~$v_0$, $v_1$, $v_2, \dots$ are assumed $2\pi$-periodic in the fast variable $\phi$).
Introducing $R(z)=\frac{7}{4}\frac{f}{f_z}-\frac{3}{2}z$, one can show that $R$ satisfies ODE (\ref{eq:Kudashev}).
This ODE was derived by Vadim Kudashev in the late 1990s, but was never published during his lifetime.
It has first appeared in~\cite{GST}, see also~\cite{SS} where a peculiar hypergeometric integral was provided, thus confirming its integrability. We also refer to~\cite{S, Dubrovin, Claeys1, Claeys2} for the universality property of system~(\ref{kdv}).

We show that equation~\eqref{eq:Kudashev} possesses the general parametric solution
\begin{gather}\label{eq:KudashevSol}
 R=\frac{\epsilon \sqrt{15}\,w}{3\sqrt{144s(s-1)w_s^2+5w^2}},\quad
z=-8\epsilon\sqrt{15}\, \frac{144s^2(s-1)w_s^3-72s(s-1)ww_s^2+\frac{5}{12}w^3}{3(144s(s-1)w_s^2+5w^2)^{3/2}},
\end{gather}
(here and in what follows $\epsilon=\pm1$), where  $w(s)$ is the general solution to the hypergeometric differential equation,
\begin{gather}\label{eq:HyperGeom}
s(1-s)w_{ss}+\left(\frac12-\frac56s\right)w_s+\frac{35}{144}w=0,
\end{gather}
corresponding to the parameter values $(\alpha, \beta, \gamma)=\left(\frac5{12}, -\frac7{12}, \frac12\right)$.
Equation~(\ref{eq:Kudashev}) also possesses a special algebraic solution that can be represented implicitly as
\begin{gather}\label{eq:KudashevSolImplicit}
20(1-3R^2)^3-27(z+14R^3-4R)^2=0.
\end{gather}

Following \cite{GST} and utilising parametric representation (\ref{eq:KudashevSol}), we calculate explicit form of the leading term $v_0$ in the asymptotic expansion (\ref{asy}). Remarkably, the requirement of  $2\pi$-periodicity of $v_0$  in the fast variable $\phi$ singles out a special (non-algebraic)  separatrix solution of the Kudashev equation that corresponds to a Frobenius solution of  the hypergeometric equation (\ref{eq:HyperGeom}).

The structure of the paper is as follows. In Section~\ref{sec:Method} we carry out symmetry reduction of a general $\mathrm{SL}(2,\mathbb R)$-invariant
third-order ODE to a first-order ODE and identify a class of Abel equations among such reductions.
Using the fact that one can construct parametric solutions of the third-order equations,
we give its analogue for the identified Abel equations.
We review the literature on integrable Abel equations and incorporate the identified class into a hierarchy of known integrable Abel equations
in Section~\ref{sec:LitReview}.
In Section~\ref{sec:Examples} we exemplify our method with some known $\mathrm{SL}(2,\mathbb R)$-invariant
third-order ODEs, in particular those satisfied by modular forms,
and elaborate on the Kudashev equation.
The phase portrait of the Kudashev equation is discussed in Section~\ref{sec:num}.
Following~\cite{GST}, in Section~\ref{sec:asy} we present the leading term~$v_0$ of the asymptotic expansion~(\ref{asy}).
In Section~\ref{sec:GenLinearisability}, we generalise the linearisability result of Section~\ref{sec:Method}.
Finally, Section~\ref{sec:Conclusion} is left for conclusions.

\section{Crux of the method}
\label{sec:Method}

Our starting point is the third-order ordinary differential equations $F(z, g, g',g'', g''')=0$ for~$g(z)$ (here prime denotes differentiation by~$z$)
that possess $\mathrm{SL}(2, \R)$-symmetry of the form
\begin{equation}
\tilde z=\frac{\alpha z+ \beta }{\gamma z+ \delta}, \quad \tilde g= (\gamma z+ \delta)^2 g+\gamma(\gamma z+ \delta), \quad
\begin{pmatrix} a & b \\ c & d\end{pmatrix} \in\mathrm{SL}(2,\mathbb R).
\label{gsym}
\end{equation}
The corresponding Lie invariance algebra  is $\mathfrak g=\langle \p_z, z\p_z-g\p_g, z^2\p_z-(2zg{+}1)\p_g\rangle$.
It turns out that the presence of symmetry (\ref{gsym}) implies linearisability of the equations under study.
The following  statement is, essentially, contained in~\cite{ClarksonOlver1996}:

{\bf Theorem 1.}  {\it A general  third-order equation  $F(z, g, g', g'', g''')=0$ possessing $\mathrm{SL}(2, \R)$-symmetry
(\ref{gsym}) can be represented in the form $F={\cal F}(I_2, I_3)=0$ where $I_2$ and $I_3$ are the basic differential invariants of  the order two and three, respectively:
$$
I_2=\frac{(g''-6gg'+4g^3)^2}{(g'-g^2)^3},   \quad
I_3=\frac{g'''-12gg''-6(g')^2+48g^2g'-24g^4}{(g'-g^2)^2}.
$$
The general solution of any such equation can be represented parametrically as
\begin{equation}
z=\frac{\tilde w}{w}, \qquad g=\frac{ww_s}{W}
\label{sub}
\end{equation}
where $w(s)$ and $\tilde w(s)$ are two linearly independent solutions of a second-order linear equation $w_{ss}+p w_s+q w=0$, 
and $W=\tilde w_sw-w_s\tilde w$ is the Wronskian of~$w$ and~$\tilde w$
(the coefficients~$p(s)$ and~$q(s)$ depend on the  equation $F=0$ and can be efficiently reconstructed, see the proof).}

\medskip

\centerline {\bf Proof:}

\medskip

\noindent  Consider a linear equation $w_{ss}+p w_s+q w=0$, take its two linearly independent solutions~$w(s)$,
$\tilde w(s)$ and introduce  parametric relations (\ref{sub}).
Using  $\mathrm ds/\mathrm dz=w^2/W, \ W_s=-pW$ and the chain rule we obtain
\begin{gather*}
g'-g^2=-q\frac{w^4}{W^2},\\
g''-6gg'+4g^3=-(q_s+2pq)\frac{w^6}{W^3},\\
g'''-12gg''-6(g')^2+48g^2g'-24g^4=-(q_{ss}+2p_sq+5q_sp+6p^2q)\frac{w^8}{W^4};
\end{gather*}
recall that prime denotes differentiation by $z$. Thus, one arrives at the relations
$$
\begin{array}{c}
I_2=\displaystyle{\frac{(g''-6gg'+4g^3)^2}{(g'-g^2)^3}=-\frac{(q_s+2pq)^2}{q^3}},   \\
\ \\
I_3=\displaystyle{\frac{g'''-12gg''-6(g')^2+48g^2g'-24g^4}{(g'-g^2)^2}=-\frac{q_{ss}+2p_sq+5q_sp+6p^2q}{q^2}}.
\end{array}
$$
To solve the equation ${\cal F}(I_2, I_3)=0$, it is therefore sufficient to find coefficients $p(s), q(s)$ such that
\begin{equation}\label{eqlin}
{\cal F}\left(-\frac{(q_s+2pq)^2}{q^3}, \ -\frac{q_{ss}+2p_sq+5q_sp+6p^2q}{q^2}\right)=0.
\end{equation}
This finishes the proof.\hfill $\square$

\bigskip

\noindent{\bf Remark.} Parametric formula (\ref{sub}) can be generalised as
\begin{equation*}
z=\frac{\tilde w}{w}, \qquad g=\frac{ww_s+rw^2}{W}
\end{equation*}
where, as in Theorem 1,  $w(s)$ and $\tilde w(s)$ are two linearly independent solutions of a second-order linear equation $w_{ss}+p w_s+q w=0$ and $W=\tilde w_sw-w_s\tilde w$ is their Wronskian. Here the coefficients~$p(s)$, $q(s)$ and $r(s)$ depend on the  equation $F$ and can be efficiently reconstructed,
see Section~\ref{sec:GenLinearisability}. Introducing an extra function $r(s)$  allows more flexibility in the construction.

\medskip

As a next step, we reduce third-order equation for $g(z)$ to a first-order equation by carrying out the symmetry reduction with respect
to the two-dimensional subalgebra $\langle\p_z,z\p_z-g\p_g\rangle$ of~$\mathfrak g$, cf.~\cite{ClarksonOlver1996,IbragimovNucci1997}.
In the new independent variable~$\omega$ and the new dependent variable~$\psi$,
\begin{equation}\label{psi}
\omega=\frac{g^2}{g'}, \quad \psi=\frac{(g')^3}{g^2(2(g')^2-gg'')},
\end{equation}
the invariants $I_2$ and $I_3$ take the form
\[
\hat I_2=\frac{(2(\omega-1)(2\omega-1)\omega\psi-1)^2}{(1-\omega)^3\omega^3\psi^2},\quad
\hat I_3=
\frac{\omega\psi_\omega-6\omega^2(\omega-1)(2\omega-1)^2\psi^3+(12\omega-7)\omega\psi^2+3\psi}{(\omega-1)^2\omega^3\psi^3},
\]
so that the reduced first-order equation for $\psi(\omega)$ can be represented as
${\cal F}(\hat I_2,\hat I_3)=0$. By construction, this equation will also be linearisable.
Using the expressions for $g, g', g''$ in terms of the linear equation $w_{ss}+p w_s+q w=0$, one can rewrite parametric formulae (\ref{psi})  as follows:
\begin{equation}\label{w}
\omega=\frac{w_s^2}{w_s^2-qw^2}, \quad \psi=\frac{(w_s^2-qw^2)^3}{w^2w_s^2(2qw_s^2+(q_s+2pq)ww_s+2q^2w^2)};
\end{equation}
here $w(s)$ is an arbitrary solution of the linear equation.

\bigskip

In what follows, we will consider a special two-parameter class of $\mathrm{SL}(2, \R)$-invariant third-order equations for $g(z)$ with a linear function ${\cal F}$, namely, $I_3+c_1I_2+c_2=0$ (only in this case the reduced equation is a first-kind Abel equation). In explicit form,
\begin{equation}\label{gc}
(g'-g^2)(g'''-12gg''-6(g')^2+48g^2g'-24g^4)+c_1(g''-6gg'+4g^3)^2+c_2(g'-g^2)^3=0.
\end{equation}
The corresponding first-order equation~$A_{c_1,c_2}$ for $\psi(\omega)$ is $\hat I_3+c_1\hat I_2+c_2=0$, explicitly,
\begin{equation}
\begin{split}\label{eq:GenIntegrableAbel}
\psi_\omega&-\omega(\omega-1)(4c_1(2\omega-1)^2-c_2\omega(\omega-1)+6(2\omega-1)^2)\psi^3\\
&+(4c_1(2\omega-1)+12\omega-7)\psi^2+\frac{3\omega-c_1-3}{\omega (\omega-1)}\, \psi=0,
\end{split}
\end{equation}
which is an Abel equation of the first kind depending on two parameters $c_1, c_2$.
Its general solution can be represented in parametric form~(\ref{w}) where~$w(s)$ is the general solution of a second-order linear equation $w_{ss}+p w_s+q w=0$
whose coefficients $p(s)$ and $q(s)$ can be recovered from the corresponding relation (\ref{eqlin}):
\begin{equation}\label{pq}
q(q_{ss}+2p_sq+5q_sp+6p^2q)+c_1(q_s+2pq)^2-c_2q^3=0.
\end{equation}
Note that we have a single constraint for the two unknown coefficients $p(s)$ and $q(s)$: this allows some flexibility in selecting a linear equation with desired analytic properties. Remarkably, in the case of (\ref{pq}), one can choose the linear equation to be hypergeometric (for regular values of~$c_1$ and~$c_2$: $c_1\neq-3/2$, $c_2\neq0$):
\begin{equation}\label{hyper}
s(1-s)w_{ss}+(\gamma-(1+\alpha+\beta)s)w_s-\alpha \beta w=0.
\end{equation}
Indeed, substituting the corresponding coefficients $p(s)=\frac{\gamma-(1+\alpha+\beta)s}{s(1-s)}, \ q(s)=-\frac{\alpha \beta}{s(1-s)}$ into~(\ref{pq}) one obtains the following relations among hypergeometric parameters $\alpha, \beta, \gamma$ and the parameters~$c_1$, $c_2$ of the Abel equation $A_{c_1,c_2}$:
\begin{gather*}
(4c_1 + 6)\gamma^2 - (4c_1 + 7)\gamma + c_1 + 2=0,\\[1ex]
(4c_1+6)(\alpha+\beta)^2 -c_2\alpha \beta=0,\\[1ex]
c_2\alpha \beta  - (8c_1+12)(\alpha+\beta) \gamma + (4c_1+5)(\alpha +\beta) + 2\gamma - 1=0.
\end{gather*}
These relations can be explicitly solved for $\gamma, \ \alpha+\beta$ and $\alpha\beta$, leading to the four cases:
\begin{gather*}
\gamma=\frac{1}{2}, \quad \alpha+\beta=0, \quad  \alpha\beta=0;\\[1ex]
\gamma=\frac{1}{2}, \quad \alpha+\beta=\frac{1}{4c_1+6}, \quad  \alpha\beta=\frac{1}{c_2(4c_1+6)};\\[1ex]
\gamma=\frac{c_1+2}{2c_1+3}, \quad \alpha+\beta=\frac{1}{2c_1+3}, \quad  \alpha\beta=\frac{2}{c_2(2c_1+3)};\\[1ex]
\gamma=\frac{c_1+2}{2c_1+3}, \quad \alpha+\beta=\frac{1}{4c_1+6}, \quad  \alpha\beta=\frac{1}{c_2(4c_1+6)}.
\end{gather*}
Thus, there can be several different hypergeometric equations linearising the same Abel equation. Note that the first case can be disregarded since it leads to the inconsistent condition $\omega=1$ in the formula~(\ref{w}).
Furthermore, hypergeometric equations in the second and the fourth cases are equivalent under the transformation $s\to 1-s,\ w\to w$.
In what follows, we will not distinguish between collections $(\alpha, \beta, \gamma)$ and
$(\beta, \alpha, \gamma)$ since they correspond to the same hypergeometric equation.

\medskip
It is important to note that besides the general solutions expressed via hypergeometric functions,  the Abel equations $A_{c_1,c_2}$  possess special {\it algebraic} solutions  given by parametric formulae~(\ref{w}) where $w$ satisfies a linear equation
$w_{ss}+pw_s+qw=0$ with constant coefficients~$p$ and~$q$. The substitution into~(\ref{pq}) gives a single relation among the parameters,
${
(6+4c_1)p^2=c_2q
}$
where without any loss of generality one can set $p=1$. Thus, the required linear equation~is
\begin{gather}\label{eq:LinearParameterisation}
w_{ss}+w_s+\frac{6+4c_1}{c_2}w=0.
\end{gather}

\section{Connection to known integrable Abel equations}\label{sec:LitReview}

After we identified the class~$A_{c_1, c_2}$ of linearisable Abel equations,
a natural question arises: how many of these integrable equations are new?
In~\cite{ChebTerrabRoche2000}, most of the known integrable Abel equations were categorised into 11~equivalence classes,
with canonical representatives and their first integrals being provided therein.
To this aim, invariance of the entire class of Abel equations $u'=f_3(x)u^3+f_2(x)u^2+f_1(x)u+f_0(x)$ under transformations
of the form $x\to\phi(x)$, $u\to\psi(x)u+\eta(x)$ was used.
It was shown in~\cite{Appell, Liouville} that $I_1=s_5^3/s_3^5$ and $I_2=s_5s_7/s_3^4$ are absolute invariants of this action.
Here~$s_3$, $s_5$ and~$s_7$ are relative invariants defined recursively~as
\begin{gather*}
s_3=f_0f_3^2+\frac13\left(\frac29f_2^3-f_1f_2f_3+f_3\frac{\mathrm d}{\mathrm dx}f_2-f_2\frac{\mathrm d}{\mathrm dx}f_3\right),\\
s_{2m+1} = f_3 \frac{\mathrm d}{\mathrm dx}s_{2m-1} - (2m - 1) s_{2m-1}\left(\frac{\mathrm d}{\mathrm dx}f_3+f_1f_3-\frac13f_2^2\right).
\end{gather*}
The procedure of relating a given Abel equation to the known integrable Abel equation based on this invariance was implemented in \textsf{Maple} shortly after that.
The second step towards the classification of integrable Abel equations was undertaken in~\cite{ChebTerrabRoche2003}
where a multi-parameter class~AIA of Abel equations containing all the above integrable Abel equations as elements was identified.
As a subclass, the class AIA contains a class~AIR of equations reducible to Riccati equations,
whose elements have an intimate connection to hypergeometric functions~\cite{ChebTerrab2004}.

It is the class~AIR that the equations in the class~$A_{c_1,c_2}$ are associated with,
which comes as no surprise since all $\mathfrak{sl}_2(\mathbb R)$-invariant third-order ODEs are reduced to Riccati equations via the symmetry-reduction procedure~\cite{RuizMuriel2017}.
In particular, the equation $A_{-\frac32,\frac2\alpha}$ is related to the classical Abel equation~$AD_\alpha$, $x^2y_x+xy^3+(x^2+\alpha)y^2=0$,
via the point transformation $\omega=x^2/(x^2+\alpha)$, $\psi=-(x^2+\alpha)^3y/(2\alpha x^3)$.
The equation $A_{c_1,0}$ is equivalent to the canonical AIR equation which appeared as Eq.~(28) in~\cite{ChebTerrab2004}
with $\alpha=1/(2a)$, $\beta=-1/(2a)$.

For regular values~$c_1$ and~$c_2$, the equation~$A_{c_1,c_2}$ possesses a first integral of the form
\begin{gather*}
I=\Psi^{-\frac{1}{a}}\frac{\hat\Psi\, {_2F_1}\left({-}\frac1{2a},\frac{\sqrt{b(b{-}8a)}{-}b}{2ab};\frac{a{-}1}{a};\Psi\right)
{+}\Psi^{-}\,{_2F_1}\left(1{-}\frac1{2a},\frac{\sqrt{b(b{-}8a)}{-}b}{2ab}{+}1;\frac{a{-}1}{a}{+}1;\Psi\right)}
{\bar\Psi\,{_2F_1}\left(\frac1{2a},\frac{\sqrt{b(b{-}8a)}{+}b}{2ab};\frac{a{+}1}{a};\Psi\right){+}
\Psi^{+}\,{_2F_1}\left(\frac1{2a}{+}1,\frac{\sqrt{b(b{-}8a)}{+}b}{2ab}{+}1;\frac{a{+}1}{a}{+}1;\Psi\right)}
\end{gather*}
where $\Psi$'s are (at most) rational functions of~$\psi$ with $\omega$-dependent coefficients, $a:=2c_1+3$, $b:=c_2$.
Recall that the Kudashev equation~(\ref{eq:Kudashev}) is known to possess a similar first integral, see~\cite{SS}.
Furthermore, such first integrals arose from different perspectives in~\cite{ChebTerrab2004} and~\cite{RuizMuriel2017}.
The former paper concerns transforming first integrals of Riccati equations to that of the associated Abel equations,
and the latter paper to that of the associated $\mathfrak {sl}_2(\mathbb R)$-invariant third-order ODEs.
Note that for the both second hypergeometric functions in the numerator and the denominator, the first three parameters are greater by one
than their counterparts on the left, and recall that this is the feature of the derivative of a hypergeometric function,
$\frac{\mathrm d}{\mathrm dz}\, {_2F_1}(p,q;r;z)={_2F_1}(p+1,q+1;r+1;z)$.
This is in accordance with the first integrals appeared in~\cite{ChebTerrab2004,RuizMuriel2017}.
Other known integrable Abel equations have similar first integrals with hypergeometric functions being replaced by other special functions
parameterised by at most one parameter, three of which are related to an element of the class~AIR.

An interesting observation is that although the above first integral contains four hypergeometric functions with various forms of
the parameters, in the special case when ${c_1=-\frac{1+6c}{2c}}$ and $c_2=-\frac{16c}{4c^2-1}$,
there are only two different forms, namely, 
\[_2F_1\left(d,d+\frac12,2d+1\right)\quad\text{and}\quad_2F_1\left(d,d+\frac12,2d\right)\]
where $d\in\{c,-c\}$.
It is known that  hypergeometric functions with such values of parameters take the explicit algebraic forms,
\begin{gather*}
_2F_1\left(d,d+\frac12,2d+1,\Psi\right)=\left(\frac12+\frac12\sqrt{1-\Psi}\right)^{-2d},\\
_2F_1\left(d,d+\frac12,2d,\Psi\right)=\frac1{\sqrt{1-\Psi}}\left(\frac12+\frac12\sqrt{1-\Psi}\right)^{1-2d}.
\end{gather*}

Although the general solution of equations in the class~$A_{c_1, c_2}$ in the form of a first integral has been known, this form is of limited use.
Solutions to a majority of integrable Abel equations in~\cite{PolyaninZaitsev2003} are given in parametric form,
which is reminiscent of the formulae \eqref{eq:KudashevSol} and \eqref{w}. Besides, the parametric formula~\eqref{eq:KudashevSol} features in the Gurevich--Pitaevskii problem,
see Section~\ref{sec:asy}.

\section{Examples}\label{sec:Examples}

In this section we discuss four examples of integrable Abel equations $A_{c_1,c_2}$ given by~(\ref{eq:GenIntegrableAbel}) that correspond to different choices of constants~$c_1$, $c_2$. The first three of them originate from the theory of modular forms, and the last example is related to the Kudashev equation.

\medskip

\noindent \textbf{Example 1: $\bf{c_1=0}$, $\bf{c_2=24}$.} In this case equation~(\ref{gc}) is the Chazy equation,
\[
g'''-12gg''+18g'^2=0,
\]
which is satisfied by the weight~2 Eisenstein series~$E_2(z)$
associated with the full modular group~$\mathrm{SL}(2,\mathbb Z)$. Setting  $g=\frac{1}{2}\frac{\triangle'}{\triangle}$ we obtain a fourth-order ODE for the modular discriminant~$\triangle$, see e.g.~\cite{Rankin1956, Maier}. The corresponding Abel equation $A_{0,24}$ is
\begin{gather}\label{eq:AbelChazy}
\psi_\omega-6\omega(\omega-1)\psi^3+(12\omega-7)\psi^2+\frac{3}{\omega}\, \psi=0.
\end{gather}
Its general solution can be represented in parametric form~(\ref{w}) where $w$ satisfies hypergeometric equation~(\ref{hyper})  with any of the following parameter values  $(\alpha, \beta, \gamma)$:
$\left(\frac{1}{12}, \frac{1}{12}, \frac{1}{2}\right)$,
$\left(\frac{1}{6}, \frac{1}{6}, \frac{2}{3}\right)$,
$\left(\frac{1}{12}, \frac{1}{12}, \frac{2}{3} \right)$.
Equation (\ref{eq:AbelChazy}) also possesses an algebraic solution given by parametric formula (\ref{w}) where $w$ satisfies the linear  equation~\eqref{eq:LinearParameterisation}, $w_{ss}+w_s+\frac{1}{4}w=0$. Taking $w(s)=a\mathrm e^{-\frac{1}{2}s}+bs\mathrm e^{-\frac{1}{2}s}$,
where without any loss of generality one can set $a=0, \ b=1$, gives
\[
\omega=-\frac{(s-2)^2}{4(s-1)},\quad
\psi=-\frac{8(s-1)^3}{s^2(s-2)^2},
\]
or in the explicit form,
\[
\psi(\omega)=\frac{(2\omega-2\Omega-1)^3}{2(\Omega-\omega+1)^2(\Omega-\omega)^2} \quad\text{where  } \Omega=\pm\sqrt{\omega(\omega-1)}.
\]
Analogous symmetry reduction of the Chazy equation was carried out in~\cite{JoshiKruskal1993}.
\medskip

\noindent \textbf{Example 2: $\bf{c_1=-1}$, $\bf{c_2=9}$.}  In this case equation~(\ref{gc}) takes the form
\begin{gather}\label{eq:E13}
g'''(g'-g^2)=(g'')^2-4g^3g''-3(g')^3+9g^2(g')^2-3g^4g'+g^6.
\end{gather}
It has appeared in the classification of integrable Euler--Lagrange equations; setting $g=\frac{f'}{f}$ one obtains a fourth-order ODE for $f$ satisfied by the Eisenstein series $E_{1,3}(z)$ \cite{FO}. The corresponding Abel equation~$A_{-1,9}$ is
\begin{gather}\label{eq:AbelE13}
\psi_\omega+\omega(\omega-1)(\omega-2)(\omega+1)\psi^3+(4\omega-3)\psi^2+\frac{3\omega-2}{\omega (\omega-1)}\, \psi=0.
\end{gather}
Its general solution can be represented in parametric form~(\ref{w}) where~$w$ satisfies hypergeometric equation (\ref{hyper})  with any of the following parameter values  $(\alpha, \beta, \gamma)$:
$\left(\frac{1}{3}, \frac{1}{6}, \frac{1}{2}\right)$,  $\left(\frac{1}{3}, \frac{2}{3}, 1\right)$, $\left(\frac{1}{3}, \frac{1}{6}, 1\right)$. Equation (\ref{eq:AbelE13}) also possesses an algebraic solution
given by parametric formula (\ref{w}) where $w$ satisfies the  linear  equation~\eqref{eq:LinearParameterisation}, $w_{ss}+w_s+\frac{2}{9}w=0$.
Taking $w(s)=a\mathrm e^{-\frac{1}{3}s}+b\mathrm e^{-\frac{2}{3}s}$,
where without any loss of generality one can set $a=1, \ b=1$, gives
\[
\omega=\frac{(2\sigma+1)^2}{2\sigma^2-1},\quad
\psi=-\frac{(2\sigma^2-1)^3}{4\sigma(2\sigma+1)^2(\sigma+1)^2}\quad \text{where}\ \sigma=\mathrm e^{-s/3},
\]
or in the explicit form,
\[
\psi(\omega)=-\frac{2(2\Omega+3\omega-2)^3}{(\Omega+\omega)^2(\Omega+2\omega-2)^2(\Omega+2)(\omega-2)} \quad
\text{where  }\Omega=\pm\sqrt{2\omega(\omega-1)}.
\]
Finally, this equation possesses the discrete symmetry $\tilde \omega=1-\omega$, $\tilde\psi=-\psi$.

\medskip

\noindent \textbf{Example 3: $\bf{c_1=-1}$, $\bf{c_2=8}$.} In this case equation~(\ref{gc}) takes the form
\[
g'''(g'-g^2)=(g'')^2-4g^3g''-2(g')^3+6g^2(g')^2.
\]
Up to a scaling factor, it has appeared in~\cite{AblowitzChakravartyHahn2006} as the equation satisfied by
the Eisenstein series~${\cal{E}}_2(z)$ of the level two congruence subgroup~$\Gamma_0(2)$ of the modular group.
The corresponding Abel equation~$A_{-1,8}$ is
\begin{gather}\label{eq:AbelAblowitz}
\psi_\omega-2\omega(\omega-1)\psi^3+(4\omega-3)\psi^2+\frac{3\omega-2}{\omega(\omega-1)}\, \psi=0.
\end{gather}
Its general solution can be represented in parametric form (\ref{w}) where $w$ satisfies hypergeometric equation (\ref{hyper})  with any of the following parameter values  $(\alpha, \beta, \gamma)$:
$\left(\frac{1}{4}, \frac{1}{4}, \frac{1}{2} \right)$,  $\left(\frac{1}{2}, \frac{1}{2}, 1\right)$, $\left(\frac{1}{4}, \frac{1}{4}, 1\right)$.
Equation \eqref{eq:AbelAblowitz} also possesses an algebraic solution,  the same as in Example~1, indeed, the corresponding linear equations~\eqref{eq:LinearParameterisation} are identical.

\noindent \textbf{Example 4: $\bf{c_1=-3}$, $\bf{c_2=24/35}$.} In this case equation~(\ref{gc}) takes the form
\[
(g'-g^2)(g'''-12gg''-6(g')^2+48g^2g'-24g^4)-3(g''-6gg'+4g^3)^2+\frac{24}{35} (g'-g^2)^3=0.
\]
We were not able to uncover its `modular' origin.
The corresponding Abel equation~$A_{-3,24/35}$ is
\begin{gather}\label{eq:AbelKudashev}
\psi_{\omega}+\frac{6}{35}\omega(\omega-1)(12\omega-5)(12\omega-7)\psi^3
-(12\omega-5)\psi^2+\frac{3}{\omega-1}\, \psi=0.
\end{gather}
Its general solution can be represented in parametric form (\ref{w}) where $w$ satisfies hypergeometric equation (\ref{hyper})
with any of the following parameter values  $(\alpha, \beta, \gamma)$:
$\left(\frac{5}{12}, -\frac{7}{12}, \frac{1}{2} \right)$,  $\left(\frac{5}{6}, -\frac{7}{6}, \frac{1}{3}\right)$, $\left(\frac{5}{12}, -\frac{7}{12}, \frac{1}{3}\right)$. With the choice $\left(\frac{5}{12}, -\frac{7}{12}, \frac{1}{2} \right)$, we explicitly have
\[
\omega=\frac{144s(1-s)w_s^2}{144s(1-s)w_s^2-35w^2}, \quad
\psi=\frac{(144s(1-s)w_s^2-35w^2)^3}{10080s(1-s)w^2w_s^2(144s(1-s)w_s^2+24sww_s+35w^2)},
\]
where $w$ solves hypergeometric equation (\ref{eq:HyperGeom}).
Equation~\eqref{eq:AbelKudashev} is related to the Kudashev equation~\eqref{eq:Kudashev}
by the  point transformation
\footnote{We slightly abuse notation here: $z$ in \eqref{eq:KudashevSol} is the independent variable of the Kudashev equation~\eqref{eq:Kudashev}, and it has nothing to do with the independent variable of the $\mathrm{SL}(2,\mathbb R)$-invariant equation for $g(z)$ at the beginning of the example.
Both the variables are denoted~$z$ in the literature and we wanted to keep the notation.}
\begin{gather}\label{eq:KudashevSol1}
 R=\epsilon\sqrt{\frac{1-\omega}{3(6\omega+1)}},\quad
z=-\frac\epsilon6\sqrt{\frac{1-\omega}{3(6\omega+1)}}\, \frac{2(1-\omega)(576\omega^2-333\omega+2)\psi+245}{(\omega-1)^2(6\omega+1)\psi},
\end{gather}
$\epsilon=\pm1$. Substituting the above expressions for  $\omega$ and $\psi$  into (\ref{eq:KudashevSol1}) we obtain parametric solution (\ref{eq:KudashevSol}) of the Kudashev equation.

Equation (\ref{eq:AbelKudashev}) also possesses an algebraic solution given by parametric formula (\ref{w}) where $w$ satisfies the  linear  equation~\eqref{eq:LinearParameterisation}, $w_{ss}+w_s-\frac{35}{4}w=0$. Taking $w(s)=a\mathrm e^{\frac{5}{2}s}+b\mathrm e^{-\frac{7}{2}s}$,
 where without any loss of generality one can set $a=1, \ b=1$, gives
\begin{equation}\label{partau}
\omega=\frac{(7\sigma-5)^2}{12(7\sigma^2+5)},\quad \psi=\frac{6(7\sigma^2+5)^3}{35\sigma(\sigma+1)^2(7\sigma-5)^2};
\end{equation}
here  $\sigma:=\mathrm e^{-6s}$. In the explicit form,
\[
\psi(\omega)=-\frac{1225}2\frac{(2\Omega+2\omega-7)^3}
{(\Omega-5\omega)^2(\Omega+7\omega-7)^2(12\Omega-35)(12\omega-7)} \quad\text{where  }
\Omega=\pm\sqrt{35\omega(1-\omega)}.
\]
Substituting~(\ref{partau}) into~(\ref{eq:KudashevSol1}) results in
\begin{equation}\label{Rzsigma}
R=\frac{\epsilon\sqrt{10}(\sigma+1)}{6\sqrt{9\sigma^2-10\sigma+5}},\quad
z=\frac{-\epsilon\sqrt{10}(\sigma-1)(\sigma^2-10\sigma+5)}{(9\sigma^2-10\sigma+5)^{3/2}},
\end{equation}
which is a parametric form of the implicit solution (\ref{eq:KudashevSolImplicit}).

\section{Phase portrait of the Kudashev equation}\label{sec:num}

Recall that solutions~$R(z)$ of the Kudashev equation~\eqref{eq:Kudashev} are given by parametric formula~\eqref{eq:KudashevSol1}
where for the general solution we choose the general solution of the associated hypergeometric equation~\eqref{eq:HyperGeom},
\begin{gather}\label{eq:HGsolution}
w(s)= a\  _2F_1\left(\frac5{12},-\frac7{12};\frac12;s\right)+b\sqrt{s}\ _2F_1\left(\frac{11}{12},-\frac{1}{12};\frac32;s\right),
\end{gather}
with arbitrary constants~$a$ and~$b$.
Depicted below is the phase portrait for the Kudashev equation.
Its apparent symmetry, associated with different values of~$\epsilon$, reflects the symmetry $z\to -z$, $R\to -R$ of the equation~(\ref{eq:Kudashev}).

\begin{figure}[ht!]
\centering
\caption{Phase portrait of the Kudashev equation}
\includegraphics[width=9cm]{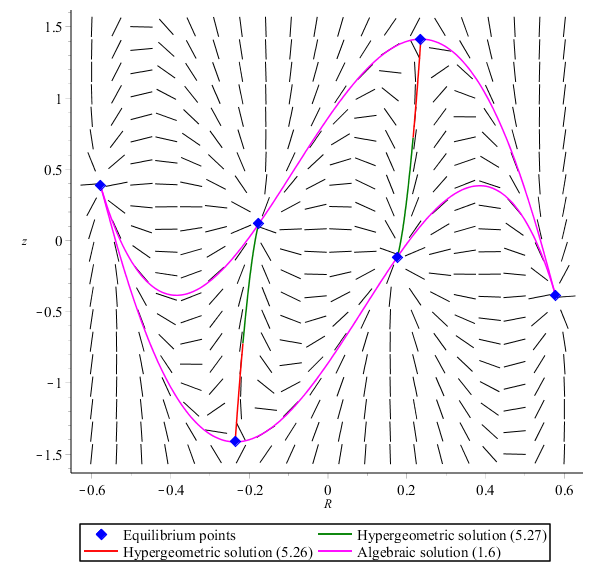}\label{fig:PP}
\end{figure}

There are six equilibrium points, where both the numerator and the denominator of the right-hand side of the equation~\eqref{eq:Kudashev} vanish
(shown in blue dots left to right in Figure~\ref{fig:PP}):
\begin{gather*}
P_1=\left(-\frac{\sqrt{3}}3,\frac{2\sqrt{3}}9\right),\quad
P_2=\left(-\frac1{3\sqrt{2}},-\sqrt{2}\right),\quad
P_3=\left(-\frac19\sqrt{\frac52},\frac2{27}\sqrt{\frac52}\right),\\[1ex]
P_4=\left(\frac19\sqrt{\frac52},-\frac2{27}\sqrt{\frac52}\right),\quad
P_5=\left(\frac1{3\sqrt{2}},\sqrt{2}\right),\quad
P_6=\left(\frac{\sqrt{3}}3,-\frac{2\sqrt{3}}9\right).
\end{gather*}
The apparent separatrix, shown in magenta in Figure~\ref{fig:PP}, is nothing else but the algebraic solution of the Kudashev equation.
It has the implicit form~(\ref{eq:KudashevSolImplicit}),
and the explicit form can be written~as
\begin{equation*}
z(R)=\frac29(3R^2-1)\sqrt{15(1-3R^2)}-2R(7R^2-2),
\end{equation*}
with different signs of the root determining the upper and the lower branches of the curve.

Another separatrix of interest, which passes through the equilibrium points~$P_4$ and~$P_5$, is not parametrised by algebraic solutions,
but by two solutions of the equation~\eqref{eq:HyperGeom},
\begin{gather}
w_1(s)=\frac{ _2F_1\left(\frac5{12},\frac{11}{12},2,\frac1s\right)}{(-s)^{\frac{5}{12}}}
=\frac{\sqrt{\pi}\ _2F_1\left(\frac5{12},{-}\frac7{12};\frac12;s\right)}{\Gamma\left(\frac{11}{12}\right)\Gamma\left(\frac{19}{12}\right)}
+\frac{2\sqrt{{-}\pi s}\ _2F_1\left(\frac{11}{12},{-}\frac{1}{12};\frac32;s\right)}{\Gamma\left(\frac{5}{12}\right)\Gamma\left(\frac{13}{12}\right)}, \label{w1}\\[1ex]
w_2(s)=\frac{\sqrt{\pi}\ _2F_1\left(\frac5{12},-\frac7{12};\frac12;s\right)}{\Gamma\left(\frac{11}{12}\right)\Gamma\left(\frac{19}{12}\right)}
-\frac{2\sqrt{-\pi s}\ _2F_1\left(\frac{11}{12},-\frac{1}{12};\frac32;s\right)}{\Gamma\left(\frac{5}{12}\right)\Gamma\left(\frac{13}{12}\right)}, \label{w2}
\end{gather}
and the value $\epsilon=1$. Note that $w_1$ is one of the Kummer solutions of~\eqref{eq:HyperGeom}.
This separatrix has a symmetric counterpart passing through the equilibrium points~$P_2$ and~$P_3$, which corresponds to the value $\epsilon=-1$.
Note that $w_1(0)=w_2(0)$, and $(R,z)|_{w=w_1}\to P_5$, $(R,z)|_{w=w_2}\to P_4$ as $s\to-\infty$.
In the subsequent figures, we always paint in magenta the lines associated with the algebraic solution of the Kudashev equation,
and in red and green the lines associated with the hypergeometric solutions of the Kudashev equations corresponding to the functions~$w_1$ and~$w_2$, respectively.

It should also be noted that integral curves lying outside the algebraic separatrix and having the endpoints $(P_4,P_6)$ and the endpoints $(P_5,P_6)$
(resp., $(P_1,P_2)$ and $(P_1,P_3)$), are also separated by a separatrix, but we were unable to find its parametrisation.

\section{Leading term of the asymptotic solution}\label{sec:asy}

Recall that the fast variable~$\phi$ in the asymptotic expansion~\eqref{asy} depends on the function~$f(z)$
that solves the first-order ODE $R(z)=\frac{7f}{4f_z}-\frac32z$, where $R(z)$ is the solution of the Kudashev equation~\eqref{eq:Kudashev}.
With the help of the parametric formula for the solution of this equation,
we can find the parametric representation of~$f(z)$.
For the general solution of the Kudashev equation, the function~$f$ takes the parametric form
\begin{gather}\label{eq:FormulaForf}
f=\frac{c|s|^{5/4}|s-1|^{5/6}|w_s|^{5/2}}{|144s(s-1)w_s^2+5w^2|^{7/4}},
\end{gather}
where $c$ is a constant to be specified later. For the algebraic solution the function $f$ has a simpler parametric form,
$f=\dfrac{c|7\sigma-5|^{5/2}}{(9\sigma^2-10\sigma+5)^{7/4}}$.

Below we follow~\cite{GST} to show that the knowledge of the coefficients~$f(z)$ and $R(z)$ leads to an explicit formula  for the leading term
$v_0(z, \phi)$ of the asymptotic expansion~(\ref{asy}); we also refer to \cite{Potemin} where an equivalent approach to this problem was developed based on the Whitham averaging procedure.  We set $v_0\equiv v$ to simplify the notation.
Substituting~(\ref{asy}) into the first equation~(\ref{kdv}), at the leading order~$t^{5/4}$ one obtains
\[
Q^2v_{\phi \phi \phi}+vv_{\phi}+Rv_{\phi}=0.
\]
Similarly, substituting~(\ref{asy}) into the second equation~(\ref{kdv}), at the leading order~$t^{3/2}$ one obtains
\[
Q^4v_{\phi \phi \phi \phi}+\frac{5}{6}Q^2\left(2vv_{\phi \phi}+v_{\phi}^2\right)+\frac{5}{18}(z-v+v^3)=0;
\]
here the coefficients $Q(z)=f_z$ and $R(z)=\frac{7}{4}\frac{f}{f_z}-\frac{3}{2}z$ are functions of $z$ only.
These two equations for~$v$ are equivalent to a single first-order equation,
\begin{equation}\label{Jaceq}
Q^2v_{\phi}^2+\frac{1}{3}v^3+Rv^2+\left(6R^2-\frac{5}{3}\right)v+5R-18R^3-\frac{5}{3}z=0.
\end{equation}
We look for a solution of~(\ref{Jaceq}) in the form
\begin{equation}\label{Jac}
v=A\, \mathrm{dn}^2\left(\frac{B}{Q}\phi, k\right)-C-R
\end{equation}
where $\mathrm {dn}(p, k)$ is the Jacobi elliptic function and the coefficients~$A$, $B$, $C$, $k$ are functions of the slow variable~$z$.
Recall that $y=\mathrm{dn}(p, k)$ satisfies the equation $y_p^2=(y^2-1)(1-k^2-y^2)$, which implies $Y_p^2= 4Y(Y-1)(1-k^2-Y)$ for $Y=\mathrm{dn}^2(p, k)$.
Substituting ansatz~(\ref{Jac}) into~(\ref{Jaceq}) we obtain four relations for the coefficients:
\begin{subequations}
\begin{gather}
A-12B^2=0,\label{eq:AsExpansionA}\\
4(2-k^2)B^2-C=0,\label{eq:AsExpansionB}\\
12(k^2-1)AB^2+3C^2+15R^2-5=0,\label{eq:AsExpansionC}\\
C^3+(15R^2-5)C+70R^3-20R+5z=0.\label{eq:AsExpansionD}
\end{gather}
\end{subequations}
One can solve the equations~\eqref{eq:AsExpansionA} and~\eqref{eq:AsExpansionB} for $A$ and $B$:
\[
A=\frac{3C}{2-k^2}, \qquad B^2=\frac{C}{4(2-k^2)};
\]
here $C$ and $k$ can be recovered from \eqref{eq:AsExpansionC} and \eqref{eq:AsExpansionD}:
\begin{equation}\label{Ck}
\begin{array}{c}
9(k^2-1)C^2+(3C^2+15R^2-5)(k^2-2)^2=0,  \\[1ex]
C^3+(15R^2-5)C+70R^3-20R+5z=0.
\end{array}
\end{equation}
For what follows, it will be convenient to rewrite equations (\ref{Ck}) in a somewhat different (equivalent) form. First of all, the elimination of $C$ from  equations (\ref{Ck}) leads to an algebraic equation for $k^2$,
\[
\frac{5(2-k^2)^2(1-2k^2)^2(1+k^2)^2}{27(k^4-k^2+1)^3}=\frac{(14R^3-4R+z)^2}{(1-3R^2)^3}.
\]
Using (\ref{eq:KudashevSol}), the right-hand side of this formula simplifies to
\[
\frac{(14R^3-4R+z)^2}{(1-3R^2)^3}=\frac{20}{27}\frac{s}{s-1},
\]
where all dependence on $w$ cancels out (note that the algebraic solution (\ref{eq:KudashevSolImplicit}) arises in the limit $s\to \infty$).
The resulting algebraic equation for $k^2$ is
\[
\frac{(2-k^2)^2(1-2k^2)^2(1+k^2)^2}{(k^4-k^2+1)^3}=\frac{4s}{s-1}, \ \  {\rm equivalently}, \ \  \frac{(2-k^2)^2(1-2k^2)^2(1+k^2)^2}{k^4(1-k^2)^2}=-27s.
\]
Secondly, solving the first equation (\ref{Ck}) for $C^2$ and substituting the result into the second equation (\ref{Ck}), one obtains an explicit formula for $C$. Ultimately, equations (\ref{Ck}) are equivalent to
\begin{equation} \label{Ck1}
s=-\frac{(2-k^2)^2(1-2k^2)^2(1+k^2)^2}{27 k^4(1-k^2)^2}, \qquad C=-3\frac{k^4-k^2+1}{(1-2k^2)(1+k^2)}\frac{14R^3-4R+z}{1-3R^2}.
\end{equation}
In what follows, we assume $k^2\in (0,1), \ s\in (-\infty, 0]$.

\medskip

The further analysis splits into two different cases  depending on whether $R(z)$ is  a generic or the algebraic solution of the Kudashev equation.

\subsection{Generic solution of the Kudashev equation}

Here we use  the generic solution (\ref{eq:KudashevSol}),
\begin{gather*}
 R=\frac{\epsilon \sqrt{15}\,w}{3\sqrt{144s(s-1)w_s^2+5w^2}},\quad
z=-8\epsilon\sqrt{15}\, \frac{144s^2(s-1)w_s^3-72s(s-1)ww_s^2+\frac{5}{12}w^3}{3(144s(s-1)w_s^2+5w^2)^{3/2}},
\end{gather*}
where $w$ satisfies hypergeometric equation~(\ref{eq:HyperGeom}).
In the Gurevich--Pitaevskii problem, it is required that the function~$v$ is $2\pi$-periodic in the fast variable $\phi$ \cite{GST},
which translates into the condition
\begin{gather}\label{eq:GPperiodicity}
\frac BQ=\frac{K(k)}\pi
\end{gather}
where~$K$ is the complete elliptic integral of the first kind.
Let us find a solution to this equation. We will use the following formulae:
\begin{itemize}

\item $B=\frac{\sqrt{C}}{2\sqrt{2-k^2}}$ where $C$ is given by the second formula  (\ref{Ck1});

\item $Q=f_z=\frac{7f}{4R+6z}$ (which follows from the definition of $R$);

\item $\displaystyle f=\frac{c(-s)^{5/4}(1-s)^{5/6}w_s^{5/2}}{(144s(s-1)w_s^2+5w^2)^{7/4}}$,\ \ $c=\operatorname{const}$,
which is a specialisation of~(\ref{eq:FormulaForf});

\item $K(k)=\frac{1}{2}\pi\, {_2F_1}(1/2,1/2,1,k^2)$.

\end{itemize}
Introducing the new independent variable $r=k^2$ (so that $r\in (0,1)$), choosing $\epsilon=-1$ and using the first formula (\ref{Ck1}) as a change of variables from $s$ to $r$, namely, $s=s(r)=-\frac{(2-r)^2(1-2r)^2(1+r)^2}{27 r^2(1-r)^2}$, we can rewrite condition (\ref{eq:GPperiodicity}) in the form
\begin{equation}\label{r}
\mu \ \frac{3r(1-r^2)(2-r)(1-2r)w_r+7(r^2-r+1)^2w}{(r(1-r))^{5/6}(r^2-r+1)}=
{_2F_1}(1/2,1/2,1,r);
\end{equation}
here the constant factor is $\mu=-\frac{128}{7c}\,  2^{5/6}\,3^{1/2}\, 5^{3/4}$ where $c$ is the same constant as in (\ref{eq:FormulaForf}). It turns out that the left-hand side of this equation indeed satisfies the hypergeometric equation with the parameters $\alpha=1/2$, $\beta=1/2$, $\gamma=1$,
whenever $w$ solves hypergeometric equation~(\ref{eq:HyperGeom}). However, we need to select a special solution of~(\ref{eq:HyperGeom}) to make~(\ref{r}) an identity. We claim that the correct choice for $w$ is the Kummer solution discussed in Section~\ref{sec:num}, namely,
\[w=(-s)^{-5/12}\ _2F_1\left(\frac5{12},\frac{11}{12},2,\frac1s\right)\] 
where we have to substitute~$s(r)$ from the first formula~(\ref{Ck1}). This can be checked directly by comparing Taylor expansions of both sides of~(\ref{r}) at $r=0$, which fixes the constant factor as $\mu=2^{-7/6}3^{-9/4}$. Comparing the two obtained expressions for~$\mu$ we recover the exact value of $c=-2^93^{11/4}5^{3/4}7^{-1}$.

\noindent{\bf Summary.} Let us bring together all the  formulae needed to calculate the first term of the asymptotic expansion
(assuming $\epsilon=-1, \ r\in(0, \, 1/2]$). We have
\[
v=\frac{3C}{2-r}\, \mathrm{dn}^2\left(\frac{K(\sqrt r)}\pi\phi, \sqrt r \right)-C-R
\]
where
\begin{gather*}
C=\frac{2\sqrt{15}\, r(1-r)(2-r)\, w_r}{\sqrt{36r^2(1-r)^2w_r^2+5(r^2-r+1)w^2}\, \sqrt{r^2-r+1}},\\
R=-\frac{\sqrt{15}\,\sqrt{r^2-r+1}\, w}{3\, \sqrt{36r^2(1-r)^2w_r^2+5(r^2-r+1)w^2}}.
\end{gather*}
Recall that the slow and the fast variables are defined as
$z=xt^{-\frac{3}{2}}$, ${\phi=t^{\frac{7}{4}}f(z)+S(z)}$;
it was shown in~\cite{GST} that the phase shift $S(z)$ for the Gurevich--Pitaevskii solution is equal to~$\pi$. Finally,  the function $f$ and the variable $z$ are defined parametrically as
\begin{gather*}
f=c\, \frac{2^{2/3}3^{5/4}(r(1-r))^{10/3}\,w_r^{5/2}}{16\left(36r^2(1-r)^2w_r^2+5(r^2-r+1)w^2\right)^{7/4}(r^2-r+1)^{3/4}},\\
z=\frac{2\sqrt{15}\, \triangle}{9\left(36r^2(1-r)^2w_r^2+5(r^2-r+1)w^2\right)^{3/2}(r^2-r+1)^{3/2}}.
\end{gather*}
where
\[
\triangle=108r^3(2-r)(1- 2r)(1 + r)(r-1)^3w_r^3-216r^2(1 - r)^2(r^2 - r + 1)^2ww_r^2+5(r^2-r+1)^3w^3.
\]
Here the function $w(r)$ is defined as
\[
w(r)=(-s)^{-5/12}\ _2F_1\left(\frac5{12},\frac{11}{12},2,\frac1s\right) \ \ {\rm where} \ \ s=-\frac{(2-r)^2(1-2r)^2(1+r)^2}{27 r^2(1-r)^2}.
\]
The function~$w$ is equivalent to the function~$w_1$, see~\eqref{w1}, since solutions of~\eqref{eq:HyperGeom} are defined up to a nonzero multiplier.
The counterpart~$w_2$ of~$w_1$, see~\eqref{w2}, also generates a solution of~\eqref{eq:GPperiodicity} with the same~$c$ and $r\in[1/2,1)$. These two solutions together form a developing bore depicted in Figure~\ref{fig:bore} below, which is the plot of the function~$v(t,x)$. Recall that $w_1(0)=w_2(0)$ and $r=1/2$ when $s=0$ (since $r\in[0,1]$), therefore the red and the green lines in Figure~\ref{fig:bore} naturally agree. For the naturality of the choice of the domain of the function~$w\in\{w_1,w_2\}$ see the remark below.

\begin{figure}[ht!]
\centering
\caption{Development of undular bore over time}
\includegraphics[width=13.5cm]{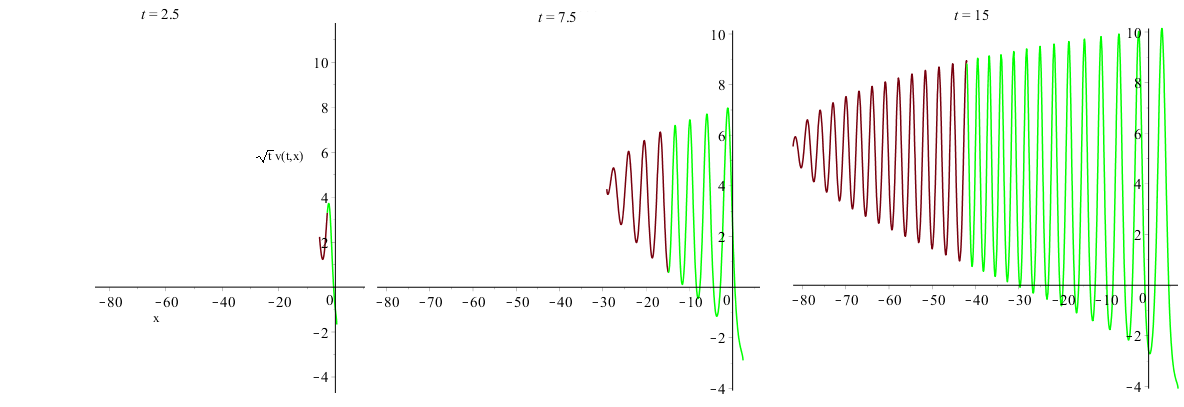}\label{fig:bore}
\caption*{The red part of the bore is parameterised by the function $w_1$ and $r\in(0,1/2]$,
the green part of the bore is parameterised by the function~$w_2$ with $r\in[1/2,1)$.}
\end{figure}

\begin{figure}[ht!]
\centering
\caption{Numerical solution of KdV equation with $u(0,x)=-x^{1/3}$ at $t=5$}
\includegraphics[width=12cm]{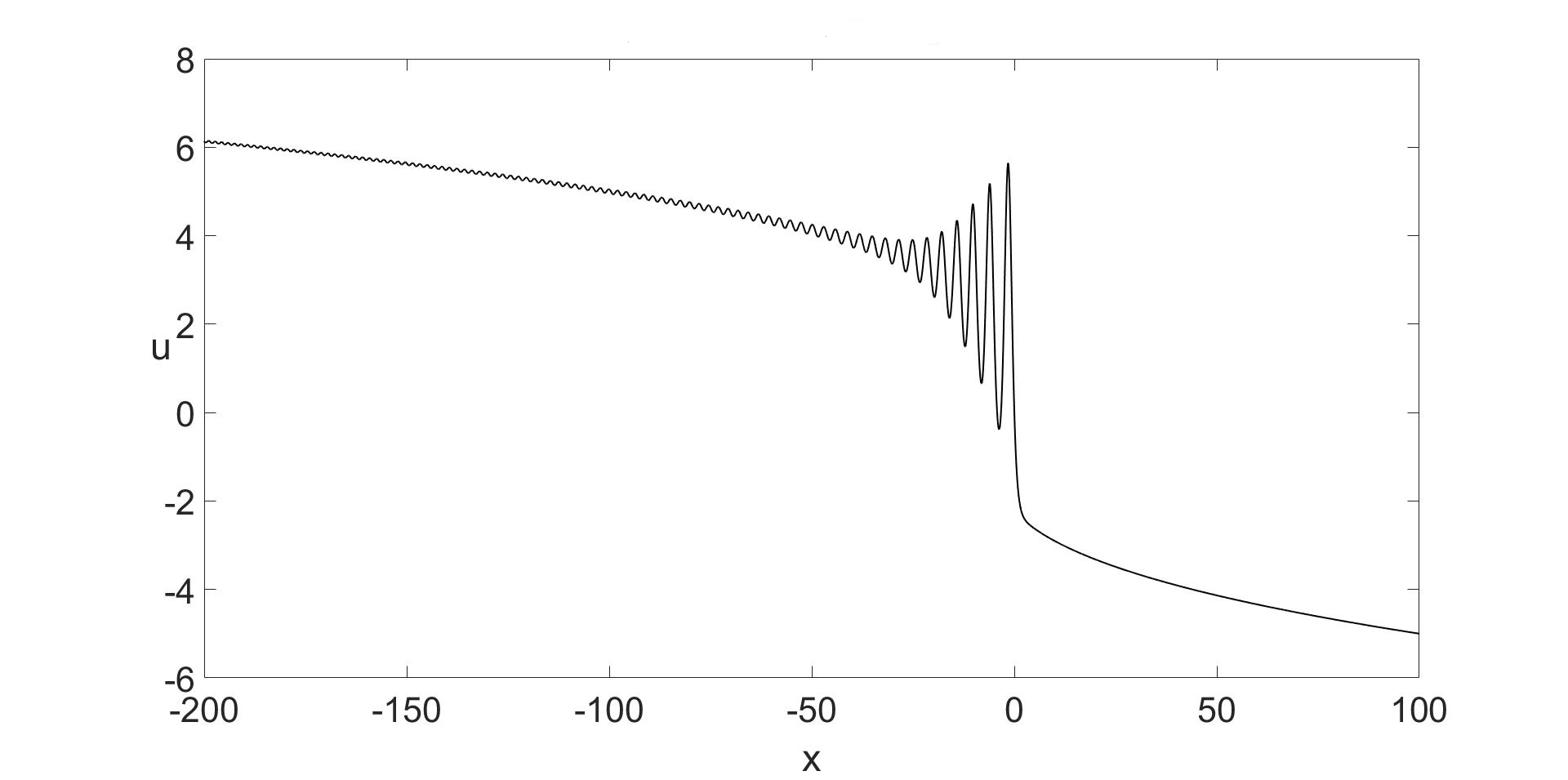}\label{fig:GP}
\end{figure}

To make the treatment of the problem comprehensive, in Figure~\ref{fig:GP}
we also include the result of numerical simulation of the KdV equation with the initial condition
$u(0,x)=-x^{1/3}$ at $t=5$ (courtesy of Curtis Hooper, see also~\cite[Fig.~5]{GP}).

\noindent {\bf Remark.} In the generic case, equation~\eqref{eq:AsExpansionD} has three distinct roots $C_1, C_2, C_3$,
which for $s\leqslant 0$ are real-valued. Using (\ref{eq:KudashevSol}) and introducing $\zeta=(\sqrt{s}+1)^{1/3}$ and  $\theta=(\sqrt{s}-1)^{1/3}$  we can represent them in the form
\begin{gather*}
C_1=
\frac{4\epsilon\sqrt{15}\, w_s(s-1)^{1/3}\sqrt{s}\left[\mathrm e^{\frac{2\pi\mathrm i}{3}}\zeta-\mathrm e^{\frac{\pi\mathrm i}{3}}\theta \right]}{\sqrt{144s(s-1)w_s^2+5w^2}},\qquad
C_2=
\frac{4\epsilon\sqrt{15}\, w_s(s-1)^{1/3}\sqrt{s}\left[\mathrm \zeta+\mathrm \theta \right]}{\sqrt{144s(s-1)w_s^2+5w^2}},\\
C_3=
\frac{4\epsilon\sqrt{15}\, w_s(s-1)^{1/3}\sqrt{s}\left[\mathrm e^{\frac{2\pi\mathrm i}{3}}\theta -\mathrm e^{\frac{\pi\mathrm i}{3}}\zeta \right]}{\sqrt{144s(s-1)w_s^2+5w^2}},
\end{gather*}
where  the corresponding values of~$k^2$ are as follows:
\begin{gather*}
k_1^2= \frac{\mathrm e^{\frac{\pi\mathrm i}3}\theta+\mathrm e^{\frac{2\pi\mathrm i}3}\zeta}{\theta-\zeta}\  \text{  and  }\
k_1^2=\frac{\theta+\mathrm e^{\frac{\pi\mathrm i}3}\zeta}{\mathrm \zeta+\mathrm e^{\frac{\pi\mathrm i}3}\theta};\qquad
k_2^2=\frac{\mathrm e^{\frac{\pi\mathrm i}3}\left(\zeta-\theta\right)}{\zeta-\mathrm e^{\frac{2\pi\mathrm i}3}\theta}\  \text{  and  }\
k_2^2=\frac{\mathrm e^{\frac{\pi\mathrm i}3}\left(\theta-\zeta\right)}{\mathrm \theta-\mathrm e^{\frac{2\pi\mathrm i}3}\zeta};\\
k_3^2= \frac{\mathrm e^{\frac{\pi\mathrm i}3}\zeta+\mathrm e^{\frac{2\pi\mathrm i}3}\theta}{\zeta-\theta}\  \text{  and  }\
k_3^2= \frac{\zeta+\mathrm e^{\frac{\pi\mathrm i}3}\theta}{\theta+\mathrm e^{\frac{\pi\mathrm i}3}\zeta}.
\end{gather*}
Recall that $k$ is the modulus of the Jacobi elliptic function, and therefore $r=k^2\in[0,1]$.
It turns out that only~$k_2^2=\frac{\mathrm e^{\frac{\pi\mathrm i}3}\left(\zeta-\theta\right)}{\zeta-\mathrm e^{\frac{2\pi\mathrm i}3}\theta}\in(0,1/2]$
and~$k_3^2= \frac{\zeta+\mathrm e^{\frac{\pi\mathrm i}3}\theta}{\theta+\mathrm e^{\frac{\pi\mathrm i}3}\zeta}\in[1/2,1)$ as $s\leqslant0$ satisfy this requirement.
The corresponding values of~$C$, $C_2$ and $C_3$, enjoy the periodicity condition~\eqref{eq:GPperiodicity} for $w=w_1$ and $w=w_2$, respectively.

\subsection{Algebraic solution of the Kudashev equation}

Here we use the algebraic solution~(\ref{Rzsigma}),
\begin{equation*}
R=\frac{\epsilon\sqrt{10}(\sigma+1)}{6\sqrt{9\sigma^2-10\sigma+5}},\quad
z=\frac{-\epsilon\sqrt{10}(\sigma-1)(\sigma^2-10\sigma+5)}{(9\sigma^2-10\sigma+5)^{3/2}}.
\end{equation*}
Note that the  implicit  equation \eqref{eq:KudashevSolImplicit}
defining the algebraic solution is equivalent to the vanishing of the discriminant  of the cubic equation~\eqref{eq:AsExpansionD}.
This means that the equation~\eqref{eq:AsExpansionD} has a multiple root.
Indeed, its roots are
\[
C_1=-\dfrac{\epsilon\sqrt{10}|7\sigma-5|}{6\sqrt{9\sigma^2-10\sigma+5}} \ \text{(of multiplicity~two)
and} \
C_2=\dfrac{\epsilon\sqrt{10}|7\sigma-5|}{3\sqrt{9\sigma^2-10\sigma+5}}.
\]
The corresponding values of~$k^2$ are~$1$ and~0, respectively. In what follows we consider the solitonic limit, $k^2=1$.
In this case, the ansatz~(\ref{Jac}) degenerates into
\begin{equation}\label{Sol}
v=A\, \mathrm{sech}^2\left(\frac{B}{Q}\phi\right)-C-R
\end{equation}
where the parameter values are as follows:
$A=3C$, $B=\frac{1}{2}\sqrt C$, $C=-\epsilon\sqrt{5/3-5R^2}$ ($\epsilon$ here determines the sign of the root).
The function~$B/Q$ is constant and equal to $-2^{5/4}3^{-3/2}5^{3/4}7^{-1}$.
Recall that $Q=f_z$ where $f=\dfrac{|7\sigma-5|^{5/2}}{(9\sigma^2-10\sigma+5)^{7/4}}$, see Section~\ref{sec:asy}.
This solution  can be interpreted as the asymptotic form of the leading soliton in the developing undular bore.

\begin{figure}[h!]
\caption{Development of a soliton solution over time}
\begin{center}
\includegraphics[width=12cm]{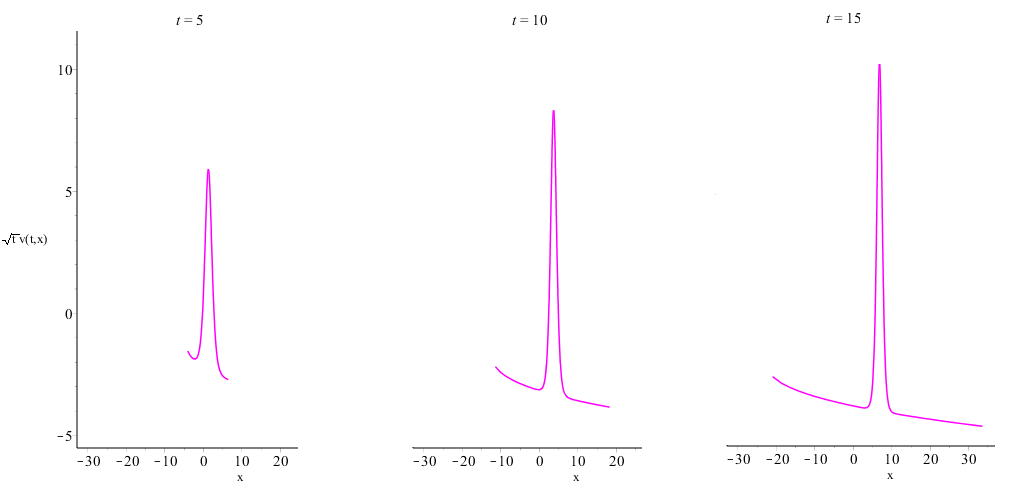}\label{fig:AS}
\end{center}
\caption*{\footnotesize{The left and the right parts of this soliton solution are parametrised differently.
One part corresponds to $\epsilon=-1$, $\sigma\in[2,1000]$, and another one corresponds to $\epsilon=1$, $\sigma\in[-1000,-2]$.
The phase shift~$S(z)$  is taken to be~0.}}
\end{figure}

It is worth to note that when $t\to\infty$, the function~$v$ tends to $V:=-C-R$, which satisfies the algebraic equation $V^3-V+z=0$, cf.~\cite{SS1}.
Recall that the behaviour of the GP special solution for $t \to -\infty$ and $x \to \pm\infty$ is principally determined from the cubic canonical equation
of the cusp catastrophe \cite{GST}
\[x - tu + u^3 = 0.\]
Changing in this equation $x=zt^{3/2}$ and $u=\sqrt{t}V$, we obtain $V^3-V+z=0$.
The situation is similar for the linear limit, $k^2=0$. In this case, the ansatz~(\ref{Jac}) degenerates into $v=\frac12C-R$ where $C=C_2$,
and $v$ again satisfies  $v^3-v+z=0$.

The determination of the phase shift~$S(z)$ requires analysis of higher-order terms in the asymptotic expansion.
As an example we take $S(z)=0$ in Figure~\ref{fig:AS} below.
The specification of the constant~$c$ in~$f$ does not seem to play any role, so it is taken~1 below.
We refer to~\cite{Grimshaw, Whitham, Dob} for the general asymptotic theory of evolution of soliton parameters and  the  phase shift problem.

\section{Generalisation of the linearisability result}\label{sec:GenLinearisability}

Here we provide the following generalisation of Theorem 1.

{\bf Theorem 2.}  {\it A general  third-order equation  $F(z, g, g', g'', g''')=0$ possessing $\mathrm{SL}(2, \R)$-symmetry~(\ref{gsym})
can be represented in the form ${\cal F}(I_2, I_3)=0$ where $I_2$ and $I_3$ are the basic differential invariants of  the order two and three, respectively:
$$
I_2=\frac{(g''-6gg'+4g^3)^2}{(g'-g^2)^3},   \quad
I_3=\frac{g'''-12gg''-6(g')^2+48g^2g'-24g^4}{(g'-g^2)^2}.
$$
The general solution of any such equation can be represented parametrically as
\begin{equation}
z=\frac{\tilde w}{w}, \qquad g=\frac{ww_s+rw^2}{W}
\label{sub1}
\end{equation}
where $w(s)$ and $\tilde w(s)$ are two linearly independent solutions of a second-order linear equation $w_{ss}+p w_s+q w=0$,
and $W=\tilde w_sw-w_s\tilde w$ is the Wronskian of~$w$ and~$\tilde w$.
Here the coefficients~$p(s)$, $q(s)$ and $r(s)$ depend on the  equation $F=0$ and can be efficiently reconstructed, see the proof below.}

\medskip

\centerline {\bf Proof:}

\medskip

\noindent  We consider a linear equation $w_{ss}+p w_s+q w=0$, take two linearly independent solutions~$w(s)$,
$\tilde w(s)$ and introduce  parametric relations~(\ref{sub1}).
Using  $\mathrm ds/\mathrm dz=w^2/W, \ W_s=-pW$ and the chain rule we obtain
$$
g'-g^2=-\tilde q\frac{w^4}{W^2},
$$
$$
g''-6gg'+4g^3=-(\tilde q_s+2\tilde p\tilde q)\frac{w^6}{W^3},
$$
$$
g'''-12gg''-6(g')^2+48g^2g'-24g^4=-(\tilde q_{ss}+2\tilde p_s\tilde q+5\tilde p\tilde q_s+6\tilde p^2\tilde q)\frac{w^8}{W^4};
$$
where $\tilde p=p-2r$, $\tilde q=q-pr+r^2-r_s$. Thus, one arrives at the relations
$$
\begin{array}{c}
I_2=\displaystyle{\frac{(g''-6gg'+4g^3)^2}{(g'-g^2)^3}=-\frac{(\tilde q_s+2\tilde p\tilde q)^2}{\tilde q^3}},   \\
\ \\
I_3=\displaystyle{\frac{g'''-12gg''-6(g')^2+48g^2g'-24g^4}{(g'-g^2)^2}=-\frac{\tilde q_{ss}+2\tilde p_s\tilde q+5\tilde p\tilde q_s+6\tilde p^2\tilde q}{\tilde q^2}}.
\end{array}
$$
To solve the equation ${\cal F}(I_2, I_3)=0$, it is therefore sufficient to find coefficients $p(s)$, $q(s)$ and~$r(s)$ such that
\begin{equation}\label{eqlin1}
{\cal F}\left(-\frac{(\tilde q_s+2\tilde p\tilde q)^2}{\tilde q^3}, \ -\frac{\tilde q_{ss}+2\tilde p_s\tilde q+5\tilde p\tilde q_s+6\tilde p^2\tilde q}{\tilde q^2}\right)=0.
\end{equation}
This finishes the proof.\hfill $\square$

Having an extra function~$r(s)$ allows one some more freedom in choosing the desired linear equation for~$w$.
For example, the general solution of the Kudashev equation can be parametrised by the associated Legendre functions~$P_{1/2}^{2/3}(s)$,
$Q_{1/2}^{2/3}(s)$ if we choose $r(s)=s^2$. In this case, it is given by the parametric formulae~\eqref{eq:KudashevSol1} where
\[
\omega=\frac{4(3(s^2-1)w_s+2sw)^2}{\Lambda(36,48,51)},\quad
\psi=\frac{\Lambda(36,48,51)^3}{280(s^2-1)w^2(3(s^2-1)w_s+2sw)^2\Lambda(-36,-36,27)},
\]
$\Lambda(\alpha,\beta,\gamma)=\alpha(s^2-1)^2w_s^2+\beta s(s^2-1)ww_s+(\gamma s^2-35)w^2$
and~$w$ is the general solution to the associated Legendre equation $(1-s^2)w_{ss}-2sw_s+\left(\frac34-\frac{4}{9(1-s^2)}\right)w=0$.
Another advantage of the parametrisation~\eqref{sub1} over the parametrisation~\eqref{sub} is that the former is invariant
under the  transformations $s\to T(s)$, $w\to S(s)w$, unlike the latter.

\section{Conclusion}\label{sec:Conclusion}

Here are a few final comments.
\begin{itemize}

\item An interesting class of exactly solvable first-order ODEs (with nonlinear dependence on the derivative) whose singular solutions can be parametrised by hypergeometric functions has appeared in~\cite{HKG} in the context of ring waves in stratified fluids (the so-called directional adjustment equations). In this connection, one should mention that algebraic separatrix solutions of the equations $A_{c_1, c_2}$ constructed in our paper can be viewed as singular solutions.

\item The algebra~$\mathfrak g=\langle\p_z,z\p_z-g\p_g,z^2\p_z-(2zg+1)\p_g\rangle$ is one of four inequivalent realisations of the Lie algebra~$\mathfrak{sl}_2(\R)$, but it is the only one that leads to Abel equations as symmetry reductions of $\mathfrak{sl}_2(\R)$-invariant third-order ODEs.
The other three realisations lead to Riccati equations~\cite{ClarksonOlver1996}.
At the same time, Abel equations are not the only equations that arise in this way.
In the case when $\mathfrak g$-invariant equation ${\cal F}(I_2,I_3)=0$ is not of the form $I_3+c_1I_2+c_2=0$,
its symmetry reduction with respect to the algebra $\langle\p_z,z\p_z-g\p_g\rangle$ is not an Abel equation,
but its solutions can still be expressed in terms of solutions of a second-order linear ODE~\eqref{eqlin},
albeit its coefficients may be hard to find explicitly.
\end{itemize}

\section{Acknowledgements}

We thank E. Cheb-Terrab, G. El, C. Hooper, R. Garifullin, K. Khusnutdinova, M. Pavlov, R.O. Popovych,  A. Shavlukov, B. Suleimanov,
S. Svirshchevskii and R. Vitolo for clarifying discussions. We also thank the referees for useful comments.
The research of SO was supported by the NSERC Postdoctoral Fellowship program.
The research of EVF was supported by a grant from the Russian Science Foundation \\
No. 21-11-00006, https://rscf.ru/project/21-11-00006/.

\end{document}